\documentclass{article}
\usepackage{graphicx}
 \usepackage{amsmath}      
 \usepackage{amssymb} 
 \usepackage{enumerate}

\newcommand{\PP}{\mathbf{P}} 
   \newcommand{\ZZ}{\mathbb{Z}} 
    
   \newcommand{\OO}{\mathcal{O}} 
    
      \newcommand{\OP}{\mathcal{O}_{\PP^3}} 
   \newcommand{\OPN}{\mathcal{O}_{\PP^N}}

   \newcommand{\II}{\mathcal{I}} 
   \newcommand{\FF}{\mathcal{F}} 
   \newcommand{\EE}{\mathcal{E}} 
   \newcommand{\GG}{\mathcal{G}} 
   \newcommand{\HH}{\mathcal{H}} 
   \newcommand{\LL}{\mathcal{L}} 
   \newcommand{\RR}{\mathcal{R}} 
   \newcommand{\SSS}{\mathcal{S}} 
    
   \newcommand{\IY}{\II_{Y}} 
   \newcommand{\gsr}{\mathrm{gsrk}} 
      \newcommand{\coker}{\mathrm{coker}} 
   \newcommand{\rk}{\mathrm{rk}}
    \newcommand{\depth}{\mathrm{depth}} 
    \newcommand{\hd}{\mathrm{hd}} 
       \newcommand{\codim}{\mathrm{codim}}

\newcommand{\im}{\mathrm{Im}\kern.3pt} 
\newcommand{\st}{\mathrm{st}} 
\newcommand{\gst}{\mathrm{gst}} 
\newcommand{\Ext}{\mathrm{Ext}} 
\newcommand{\Hom}{\mathrm{Hom}} 
\newcommand{\freccia}{\longrightarrow} 
\newcommand{\implica}{\Longrightarrow}

  \newtheorem{lemma}{Lemma}[section]
    \newtheorem{example}[lemma]{Example}
     \newtheorem{theorem}[lemma]{Theorem}
     \newtheorem{corollary}[lemma]{Corollary}
     \newtheorem{proposition}[lemma]{Proposition}

     \newtheorem{remark}[lemma]{Remark}
     
     \newtheorem{definition}[lemma]{Definition}

\newenvironment{proof}{ \noindent {\bf Proof:} }{$\diamond $ 
\vspace{0.5cm}}

\DeclareMathAlphabet{\mathpzc}{OT1}{pzc}{m}{it} 

\begin{document}

  \title{Positivity of Chern classes for Reflexive Sheaves on $\PP^n$}

      \date{}

\pagestyle{myheadings}
\markboth{ }{ }
     \author{ Cristina Bertone-Margherita Roggero \\
            Dipartimento di Matematica dell'Universit\`{a} \\
            Via Carlo Alberto 10 \\
            10123 Torino, Italy \\  
            {\small margherita.roggero@unito.it %
}          \\{ \small cristinabertone@hotmail.it}
 }

\maketitle

\footnotetext{Mathematics Subject Classification 2000: 14F05, 14F17, 14Jxx
\\ Keywords:  rank 2 bundles, subcanonical surfaces.}

\begin{abstract} 
It is well known that the Chern classes $c_i$ of a rank $n$ vector bundle  on $\PP^N$, generated by global sections, are non-negative if $i\leq n$ and vanish otherwise. This paper deals with the following question: does the above  result hold for the wider class of reflexive sheaves? We show that the Chern numbers $c_i$ with $i\geq 4$ can be arbitrarily negative for reflexive sheaves of any rank;  on the contrary for $i\leq 3$ we show positivity of the $c_i$ with weaker hypothesis. We obtain lower bounds  for $c_1$, $c_2$ and $c_3$ for every reflexive sheaf $\FF$ which is generated by  $H^0\FF$ on some non-empty open subset and completely classify sheaves for which either of them reach the minimum allowed, or some value close to it.
\end{abstract}

\section{Introduction}

In this paper we investigate some general conditions that ensure the 
positivity of Chern classes for reflexive sheaves of any rank on the projective space $\PP^N$.
There are some classical  results about  vector bundles:

\medskip

\noindent \textit{If a rank $n$ vector bundle $\FF$ on $\PP^N$ is generated by global 
sections, then  its Chern classes $c_i$ are  non-negative  if $i\leq  n$, 
while the following ones vanish} (see \cite{fulton}, Example 12.1.7).

\medskip

We would like to weaken both   hypotheses, considering the wider class 
of  reflexive sheaves (instead of vector bundles)  generated by  global 
sections on some (non-empty) open subset. 

In this new context, the 
situation is immediately more complicated. First of all, a rank $n$ reflexive 
sheaf has in general non-zero Chern classes $c_i$ also for $i>n $.
Moreover,  it is not difficult to obtain for every pair $(n,i)$ 
(both $\geq 4$)   examples of rank $n$ reflexive sheaves on $\PP^N$, generated 
by global sections,  having  negative $c_i $ (see Example \ref{negativo}).

So,
 we can not expect to control the 
positivity of the $i$-th Chern class for every    reflexive sheaf when  
$i\geq 4$, even if $i$ is lower than the rank. 

The  different behavior  of a general reflexive sheaf $\FF$ with respect to a vector bundle is clearly due to the presence of its  \lq\lq singular locus\rq\rq\ $\SSS$, that is the set of points where $\FF$ is not  locally free;  $\SSS$  is a closed subset of
  codimension  $\geq 3$,  so that if $i> 3$  the $i$-th Chern 
class $c_i(\FF)$, which is given by a    cycle  of codimension $i$,  can have components contained in $\SSS$. So  we cannot expect that $c_i(\FF)$ is necessarily positive, even if  the locally free sheaf  $\FF_U$, the restriction of $\FF$ to  $U=\PP^n\setminus \SSS$, is   generated by global sections. 

We might  think to apply the same argument 
   also to the third Chern class $c_3$ and even to the  lower ones  $c_1$ and $c_2$ in case the
     reflexive sheaf 
(or even the bundle)   is not globally generated on some closed subset of small codimension.   Thus,   it is a little surprising 
to discover that, on the contrary,   $c_1$,  $c_2$ and $c_3$ are positive under the above weaker conditions. 
 
 In fact, in \S 3 we obtain the above quoted positivity results for the first and second Chern  classes of a rank $n$ reflexive sheaf  $\FF$, generated by global sections on a non-empty
 open
 subset of $\PP^N$, as a consequence of  more general inequalities involving $c_1$ and $c_2$. More precisely,:
 
 \noindent \textbf{Theorem A} \textit{If  $\FF$ is not a direct sum of line bundles  and it has a proper subsheaf  isomorphic to  $\oplus_{i=1}^n \OPN (\alpha_i)$ where  
 $\alpha_1\geq   \dots \geq \alpha_n \geq 0$, then: 
$$
c_1\geq   \sum  \alpha_i +1  \hbox{ \ \ \ and \ \ \ \ }  c_2(\FF)\geq \sum_{ i< j}\alpha_i\alpha_j +\sum_{i\neq 2} \alpha_i+1.
$$
Moreover
$$
c_1\geq   \sum \alpha_i +2  \hbox{ \ \ \ and \ \ \ \ } c_2(\FF)\geq \sum_{ i<j}\alpha_i\alpha_j +2\sum_{ i\neq 2} \alpha_i+2
$$
unless $\FF$ has a short free resolution of the type:}
$$
0\freccia \OPN(\beta-1) \freccia \oplus_{i=1}^n \OPN (\alpha_i)\oplus \OPN(\beta) \freccia 
\FF \freccia 0 
$$ 
(see Corollary \ref{c1+1}  and Corollary  \ref{boundperc2}).

 In \S 4 we obtain similar results on $c_1$ and $c_2$ using a slightly different set of hypotheses, also involving the general splitting type of $\FF$.
 
 Finally in \S 5 we obtain similar results about the    third Chern class $c_3$ of a rank $n$ reflexive sheaf on $\PP^N$:
 
\medskip

 \noindent \textbf{Theorem B} \textit{if $\FF$ is generated by global sections outside a closed subset of codimension $\geq 3$, then $c_3(\FF)\geq 0$ 
and equality $c_3(\FF)=0 $ can hold only
 if either $N=3$ and $\FF$ is a vector bundle or $N\geq 4$ and  $\FF_H$ is a vector bundle
for every   general linear subspace $H\cong \PP^3$   in $\PP^N$}
\medskip

\noindent (see Theorem \ref{positivc3}). Under some additional condition on the homological dimension of $\FF$,  $c_3$ can vanish only if $\FF$ is a bundle:

\medskip
 
 \noindent \textbf{Theorem C} \textit{ if  $\hd(\FF)\leq 1$ and $\hd(\FF^\vee)\leq 1$, then $
c_3(\FF)=0 $ only if $\FF $ is a vector bundle having  a direct summand $ \OPN^r$, for some $r\geq n-2-h^1\FF(-c_1)$}

\medskip

\noindent (see Corollary \ref{semprec3}). This  extends  to sheaves of any rank on projective spaces of any dimension a well known property for rank $2$ reflexive sheaves on $\PP^3$   (see \cite{HSRS}, Proposition 2.6).

\section{Notation and preliminary results}

In this paragraph we introduce some general facts on reflexive sheaves and Chern classes that we will use in the paper and especially we study the singular locus of a reflexive sheaf $\FF$ and its maximal free subsheaves  (see Definition \ref{typeglobalsections}) that will be the main tool in the proofs. 

In what follows, we consider an algebraically closed field $k$ of characteristic $0$. Actually, the results of \S 3 hold more generally over a ground field of any characteristic, while we need the characteristic $0$ in \S 4 and \S 5 in order to use Grothendieck's Theorem for vector bundles on a line and Generic Smoothness.

   $\PP^N$ is the projective space of dimension $N$ over $k$. As usual, if $\FF$ is a coherent sheaf on $\PP^N$, we will denote by $h^i(\FF)$  the dimension of the $i$-th cohomology module $H^i(\FF)$ as a $k$-vector space and by $H^i_*\FF$ the direct sum $\oplus_{n\in \ZZ} H^i\FF(n)$; in particular $H^0_*\OPN=k[X_0,\dots,X_n]$ and for $\FF$ coherent sheaf, $H^0_*\FF$ has a natural structure of $H^0_*\OPN$-module; if $Y$ is a subvariety (that is a closed subscheme) in $\PP^N$, we will denote by $\deg_2(Y)$ the degree of the codimension $2$ (may be reducible or not reduced) component of $Y$.

\medskip

We recall some basic properties of Chern classes and reflexive sheaves. 
\begin{enumerate} 
\item For every     coherent sheaf $\FF$  on $\PP^N$, we denote by $c_i(\FF)$ or 
simply $c_i$ ($i=1, \dots N$) its Chern classes  that we think  as integers 
and by 
\[ 
C_t(\FF)=1+c_1(\FF)t+\dots+c_{N-1}(\FF)t^{N-1}+c_N(\FF)t^N 
\] its Chern polynomial. 
If 
$0\freccia \FF'\freccia \FF \freccia \FF''\freccia 0$ is an exact sequence, then 
$ 
C_t(\FF)=C_t(\FF')C_t(\FF'') 
$ in $  \ZZ[t]/(t^{N+1})$. 

If $c_i$ are the Chern classes  of a rank $r$ coherent sheaf $\FF$,  the Chern classes 
of $\FF (l)$ are given by: 
\begin{equation}\label{cherncontwist} 
c_i(\FF(l))=c_i+(r-i+1)lc_{i-1}+{r-i+2\choose 2}l^2c_{i-2}+\cdots+{r\choose i}l^i. 
\end{equation}

\item   We say that a coherent sheaf $\FF$ on $\PP^N$  is \emph{reflexive} if 
the canonical morphism $\FF\freccia\FF^{\vee\vee}$ is an isomorphism, where $\FF^{\vee}$ is the dual sheaf, that is 
$\FF^{\vee}=\mathcal{H}om(\FF, \OPN )$.  We refer to \cite{HSRS} and to 
\cite{okted} for general facts about reflexives sheaves, especially about rank 
$2$ reflexive sheaves on $\PP^3$ and $\PP^4$. We only  recall some of them that 
we will use more often. The dual of every sheaf is reflexive and, for every 
integer $l$,  $\OPN (l)$ is the  the only rank 1 reflexive sheaf on $\PP^N$ 
with $c_1=l$. 

A reflexive sheaf $\FF$ is locally free except at most on a closed subset 
$\mathcal S(\FF)$ of codimension $\geq 3$,   its \textbf{singular locus}.  Then, 
reflexive sheaves on $\PP^1$ and $\PP^2$ are in fact vector bundles, while on $\PP^N$ if $N\geq 3$ there 
are reflexive sheaves which are not vector bundles:  for every irreducible,  codimension 2, subvariety $Y$ in $\PP^N$ and a general section of $\omega_Y(t)$ ($t\gg 0$) we can construct a non-locally free reflexive sheaf   of rank $2$ as an extension:
$$0 \rightarrow  \mathcal{O}_{\mathbb{P}^N} \rightarrow \mathcal{F} \rightarrow \mathcal{I}_Y(b) \rightarrow 0$$ where $b$ depends on $Y$ and $t$ (see  \cite{HSRS}, Theorem 4.1). 

In the following we will use several times the following general facts:
\end{enumerate} 

\begin{lemma}\label{generale} Let  $\FF$ be a torsion free (respectively: reflexive, locally free) sheaf on $\PP^N$.   If $H$ is a  general linear subspace  of dimension $r$ in $\PP^N$, ($1\leq r \leq N-1$), then : 
\begin{description} 
\item[(i)] $H$ is \lq\lq regular\rq\rq\ with respect to $\FF$ that is $\mathrm{Tor}^1(\FF,\OO_H)=0$
\item[(ii)]    $\FF_H$ is a torsion free (reflexive, locally free) sheaf on $H$; 
\item[(iii)] the dual of $\FF_H$ as a sheaf on $H$ is isomorphic to   $(\FF^{\vee})_H$;
\item[(iv)] if $\FF$ is reflexive, the singular locus $\mathcal S(\FF_H)$ as a sheaf on $H$ is precisely $\mathcal S(\FF)\cap H$.
\item[(v)] for every $i\leq r$, $c_i(\FF_H)=c_i(\FF)$ (where $\FF_H$ is considered as a sheaf on $H$).
\end{description}
\end{lemma}
\begin{proof} 
It is sufficient to prove the results for   a general hyperplane $H$ and use induction on $N$. 

    \textit{(i)} and  \textit{(ii)} are special cases of  more general statements proved in \cite{HL} (see  Lemma 1.1.12 and Corollary 1.1.14 \textit{iii)}). In the present context $\FF$ has \lq\lq codimension\rq\rq\  $c=0$, so that \lq\lq $\FF$ pure (reflexive) of codimension 0\rq\rq\ is equivalent to $\FF$ torsion free (reflexive), and \lq\lq $\FF_H$  pure (reflexive) of codimension 1 as a sheaf on $\PP^N$ \rq\rq\ is equivalent to $\FF_H$ torsion free (reflexive) as a sheaf on $H$. 

For \textit{(iii)}, that is the isomorphism $(\FF^{\vee})_H\cong (\FF_H)^{\vee}$, see  again \cite{HL}, the remark after Definition 1.1.7.

\textit{(iv)} As in \textit{ii)} we can see that  $\mathcal S(\FF_H)\subseteq \mathcal S(\FF)\cap H$, because,  for a general $H$, $\FF_H$ is locally free (as a sheaf on $H$) where $\FF$ is (as a sheaf on $\PP^N$). Moreover $\PP^N$ is a regular variety and $H$ is  regular with respect to $\FF$;  then for every point $x$ contained in the hyperplane $H$, $\dim (\mathcal{O}_{\PP^N\!  ,x})- \depth(\FF_x)$ as  $\mathcal{O}_{\PP^N\!  ,x}$-modules, coincides with $\dim (\mathcal{O}_{H,x})-\depth(\FF_{H,x})$ as   $\mathcal{O}_{H,x}$-modules, so that  $x$ is a singular point for $\FF$ if and only if it is for $\FF_H$.

\textit{(v)}  Fix a free resolution of $\FF$; for a general $H$ it restricts to a free resolution of $\FF_H$ on $H$. Now it is sufficient to observe that the  equality   holds for all free sheaves and use  the multiplicativity of Chern polynomials.
 
\end{proof}

\begin{lemma}\label{smooth} Assume that the ground field $k$ has characteristic 0.
Let $\FF$ be a rank $r$ reflexive sheaf on $\PP^N$ generated by global sections outside  a closed subset $Z$ of
 codimension  $\geq 2$. Then 
$n-1$ general global sections degenerate on a closed subset $Y$ of codimension $\geq 
2$, generically smooth outside $Z$ and ${\mathcal S}$. 
\end{lemma} 

\begin{proof} 
On the open subset $U=\PP^N-( Z\cup {\mathcal S})$  of $\PP^N$  the  restriction map $H^0\FF \rightarrow H^0\FF_U$ is in fact a bijection 
(\cite{HSRS}, Proposition 1.6) and so the bundle $\FF_U$ on $U$  is generated by its global sections too (so that   $h^0\FF=h^0\FF_U\geq r$). 

 Take $n=r-1$ general global sections $s_1, \dots, s_{r-1}$ and consider their  degeneracy locus   $Y'=\{ x \in U \hbox{ s.t. } \dim Span(s_1(x), \dots, s_{r-1}(x))\leq r-2\} $ on $U$. If  both $m=r-2$ and $m=r-3$ satisfy the  inequality   $max\{ 0, 2r-1-h^0\FF_U\}\leq m \leq r-1$, we can  apply Remark 6 of \cite{KL}  to the bundle $\FF_U$ and to the vector space $V=H^0\FF_U$ and conclude that either $Y'$ is  empty or  it has pure codimension $2$ and it is smooth outside the subset   $Y''=\{ x \in U $ s.t. $ \dim Span(s_1(x), \dots, s_{r-1}(x))\leq r-3\} $ of codimension  $\geq 3$.

If  $h^0\FF_U=r$,  the degeneracy locus on $U$ of $r-1$ general sections is empty. If $h^0\FF_U=r+1$,  we can deduce from \cite{KL}  that  either $Y'$ is empty or  $\codim(Y')=2$ and $Y'$ is smooth outside $Y''$. In order to prove that $\codim(Y'')\geq 3$ we can   apply \cite{fulton},  Example 14.3.2 (d) to $\FF_U$ with      $p=h=3$ and $\lambda'=(3,0 \dots,0)$ (observe  that the first $r-1$ of $2r$ general sections are general too).

Finally, if $Y$ is the degeneracy locus on $\PP^N$ of $s_1, \dots, s_{r-1}$, then $Y\cap U= Y'$, so that   $\codim (Y)\geq 2$ and it is smooth outside $Z$,   ${\mathcal S}$  and $Y''$.
 
\end{proof} 

\begin{enumerate}

\item[3.] Assume that the ground field $k$ has characteristic 0. For any rank $n$ reflexive sheaf $\FF$,  we denote by $st(\FF)$ the  {\bf 
splitting type } of $\FF$ that is the sequence of integers $(b_1, \dots , b_n)$, (we will assume $b_1\geq \dots\geq b_n$)  such that for every 
general line $L$ in $\PP^N$, $\FF_L = \oplus_{i=1}^n \OO_L (b_i)$; recall that 
$c_1(\FF)=b_1 + \cdots +b_n$.

\end{enumerate}

\begin{lemma}\label{spezzecodim2} 
Let $\FF$ be a reflexive sheaf on $\PP^N$ with Chern classes $c_i$ and let 
$Y$ be a subvariety of $\PP^N$ of codimension $\geq 2$. If there is an exact 
sequence 
\begin{equation}\label{1seqesatta} 
0\freccia \GG\freccia \FF\freccia \IY(q)\freccia 0 
\end{equation} 
then $\GG$ is  reflexive and : 
\begin{equation}\label{uguaglianze}  c_1(\FF)  =c_1(\GG)+q  \hspace{1cm}   
c_2(\FF)  =c_2(\GG)+\deg_2(Y)+c_1(\GG)q.\end{equation}    
Furthermore, if $\GG$ is a vector bundle, then $Y$ is empty or it has pure 
codimension 2. 
\end{lemma} 

\begin{proof} The equalities on $c_1$ and $c_2$ can  be easily computed by 
(\ref{1seqesatta}) using 1.\  in \S 2: in fact $c_1(\FF)=c_1(\GG)+c_1(\IY(q))$ and 
$c_2(\FF)=c_2(\GG)+c_2(\IY (q))+c_1(\GG)c_1(\IY(q))$ where 
$c_1(\IY(q))=q$ and $c_2(\IY(q))=\deg_2(Y)$. 

In order to show that $\GG$ is reflexive it is sufficient to observe that $\IY$ 
is torsion-free and use \cite{HSRS} Corollary 1.5.

Finally, suppose that $\GG$ is a vector bundle and let $p$ be a codimension $\geq 3$ 
point in $\PP^N$. Then $\depth ({\GG}_p )=\codim(p)\geq 3$ and $\depth ({\FF}_p)\geq 2$, because  $\GG$ is locally free and $\FF$ is reflexive (see \cite{HSRS} 
Proposition 1.3).
 Lo\-ca\-li\-zing   the exact sequence (\ref{1seqesatta}) at the point 
$p$ and using the characterization of depth in terms of 
 non-vanishing of ${\rm {Ext}}^{i}$ (\cite{Mat} 
  Theorem 28), we  see that 
$\depth (\II_{Yp}) =\depth (\FF_p)\geq 2$, so that $\depth (\OO_{Yp})\geq 1$  so that   $p$ is
 not an
 associated prime to $Y$. 
 
\end{proof}

\begin{definition} 
Let $\HH$ be a coherent sheaf on $\PP ^N$. \\ 
We say that $\HH$ has $m$ {\bf independent global sections} if there is an 
injective map: 
\[ 
\phi :   \OPN^{m} \rightarrow \HH. 
\] 
We will denote by  $\gsr (\HH)$ the   maximum $m$ for  which   $\HH$ has $m$ 
independent global sections and call  this number the \textbf{ global section rank} of $\HH$. 
\end{definition} 

If $m=\gsr(\HH)$,   there are global sections $s_1,\dots s_m  $  in $ H^0\HH$ 
which are linearly independent  on $H^0_*\OPN=k[X_0,\dots,X_N]$. Especially, 
if $m=\gsr(\HH)=\rk(\HH)$, any set of $m$  independent   global sections 
generate $\HH$ outside the hypersurface zero locus of $s_1 \wedge \dots 
\wedge s_m$.  Of course, $\gsr (\HH ) \leq \rk(\HH)$ and equality holds if, but 
not \lq\lq only if \rq\rq , $\HH$ is generated by global sections.

\begin{example} 
 Let $X$ be a complete intersection $(1,d)$ ($d\geq 2$) in $\PP^3$  and let 
$\HH$ be the rank $2$ reflexive sheaf defined as a (suitable) extension: 
\[ 0\freccia \OP  \freccia \HH \freccia \II_X(q)\freccia 0 \] for some $q\ll 0$. 
Then $\gsr (\II _X (1))=1=\rk (\II _X (1))$ and $\gsr (\HH (1-q))=2=\rk(\HH 
(1-q))$, but they are not generated by their global sections. 
\end{example}

As usual, we will denote by $\hookrightarrow$ any injective map. 

\begin{lemma} 
Let $\HH$ be a torsion free sheaf on $\PP ^N$ and consider   the subsheaf 
$\EE$   generated by  $H^0 \HH   $.   For every integer $n$, the following are equivalent: 
\begin{enumerate} 
\item $\gsr(\HH)\geq n$; 
\item   $\rk(\EE)\geq n$; 
\item $\rk(\EE_L)\geq n$ as an  $\OO _L$-module,  for  a general line $L$  in $\PP 
^N$; 
\item there is a map $\phi \colon \OPN^n \rightarrow \HH$ such that $\phi_L$ 
is injective for  a general line $L$  in $\PP 
^N$. 
\end{enumerate} 

Then $\gsr(\HH)=\rk(\HH)$ if and only if $H^0\HH$ generates $\HH$ in every point of a suitable open subset $U$ of $\PP^N$
\end{lemma} 

\begin{proof}\ 

\begin{description} 
\item[$1. \implica 2. $] Let  $\mathcal B$ be the image of an injective 
map $ 
\phi: \OPN^n\hookrightarrow \HH $. Then  $\rk(\EE) \geq  \rk(\mathcal{B} )= n$. 
\item[$2. \Longleftrightarrow 3.$] For a general line $L$,  the rank of 
$\EE_L$ as an $\OO_L$-module coincides with  the rank of  $\EE$ as an 
$\OPN$-module. 
\item[$3. \implica 4.$] Fix a  base $s_1, \dots, s_r$ for  $H^0\HH$ as 
a $k$-vector space   and consider the corresponding surjective map $\psi 
\colon  \OPN^r \rightarrow \EE$.   As $\EE$ is torsion free and $L$ is 
general, then 
$\psi_L \colon  {\OO}_L^r \rightarrow \EE_L \simeq \oplus_{i=1}^m 
\OO_L(a_i)  $ is surjective too and  $m=\rk(\EE_L) \geq n$ by hypothesis. 
We get a map $\phi$ as required, for every choice of $n$ global section of $\HH$ 
such that their restriction 
to $L$ are independent. 

\item[$4. \implica 1.$]  The map $\phi$  induces an exact 
sequence: 
\[ 0 \rightarrow \mathcal{R} \rightarrow \OPN^n \rightarrow \EE' \rightarrow 
0\] where both $\EE'=Im(\phi)\subseteq \HH$  and $\RR=Ker(\phi)\subseteq \OPN^n $ are torsion free. Then for a general line $L$, we have  $\mathcal{T}or^1(\EE', 
\OO_L)=0$ and $\rk(\mathcal{R})=\rk(\mathcal{R}_L)=0$ so that  $\phi$ is injective and  $\gsr(\HH)\geq n$. 
\end{description} 

Finally if $n=\rk(\HH)=\gsr(\HH)$,  then $\HH$ and its subsheaf $\EE$ have the same rank and so they coincide on some open, non-empty, subset $U$. On the converse, if $H^0\HH$ generates $\HH$ on $U$, we can take any point $x\in U$ such that $\HH_x$ is free and $n=\rk(\HH)$  sections in $H^0\HH$ defining  a map $\phi :   \OPN^{n} \rightarrow \HH$ which is injective in $x$; then $Ker(\phi)$ is a torsion subsheaf of $\OPN^{n}$, that is $Ker(\phi)=0$.
 
\end{proof}

\begin{definition}  \label{typeglobalsections} 
Let $\FF$ be a rank $n$ torsion free sheaf on $\PP ^N$. The \textbf{ 
global section type} $\gst (\FF)$ is the sequence of 
integers $(a_1, \dots, a_n)$ such that $a_1\geq \dots\geq a_n$ and for every 
$i=1, \dots , n$ \[\gsr(\FF (-a_i -1)) <i 
\hbox{\hspace{.4cm} and \hspace{.4cm} } \gsr(\FF (-a_i )) \geq i.\] 
 
\end{definition} 

Note that  $a_i\geq 0$ if and only if $\gsr (\FF) \geq i$; moreover $\gst(\FF 
(l))=(a_1+l, \dots, a_n+l)$.

\begin{remark}\label{mdp}  The global section rank of $\FF$ is strictly related to some \lq\lq maximal free subsheaves\rq\rq\  of $\FF$, in the following sense:

  $\gst (\FF)=(a_1, \dots, a_n)$ if and only if there is
 a (not unique) injective map:  \[\phi \colon \oplus_{i 
=1}^n\OPN (a_i) \hookrightarrow \FF\] and for every injective map
 $f \colon \oplus_{i 
=1}^r\OPN (\alpha_i)\hookrightarrow \FF$  ($\alpha_1 \geq \dots \geq 
\alpha_r$), we have  $r\leq n$ and $ \alpha _i \leq a_i $  for every  $i=1, \dots , r$.

Maximal free subsheaves of a  sheaf $\FF$ are studied in \cite{MDP}, Chapter IV. However,  the two notions of maximality are slightly different  and the present one  is in one sense     weaker and in another sense stronger than  that given by \cite{MDP}. It is weaker because the sheaf  $\phi(\oplus_{i 
=1}^n\OPN (a_i))$ and its direct summand $\phi(\oplus_{i 
=1}^r\OPN (a_i))$ ($r<n$) are maximal among the free subsheaves of $\FF$ of the same rank and not among all the subsheaves of $\FF$ of the same rank as in \cite{MDP}. It is stronger because it is stable under isomorphisms, while in \cite{MDP} it is not. The following example illustrates a few differences between the two notions of maximality. 
\end{remark} 

\begin{example} \label{maximality} Let $\FF$ be a non-totally split reflexive sheaf with $\gst(\FF)=(a_1,\dots , a_n)$.
\begin{description}
\item[1)] $\oplus_{i =1}^n\OPN (a_i)$ cannot be a maximal free subsheaf of $\FF$ in the sense of \cite{MDP} because the only maximal subsheaf of $\FF$ of  rank   $n=\rk(\FF)$ is obviously $\FF$ itself.

\item[2)] every isomorphic image of $\OPN(a_1)$ is a maximal subsheaf of $\FF$ in the sense of \cite{MDP}. On the other hand, if $a_1<a_2$, the subsheaf $\phi(\OPN(a_2))$ is  maximal, but there are subsheaves of $\FF$ isomorphic to $\OPN(a_2)$ that are not, namely  $f\phi(\OPN(a_1))$, where $f$ is any form of degree $a_2 -a_1$.

\item[3)] $\oplus_{i =1}^{n-1}\OPN (a_i)$ is not always a maximal subsheaf of $\FF$ of rank $n-1$. For instance, let $\EE$ be a non-split bundle of rank $n-1$ with $\gst(\EE)=(a_1,\dots,a_{n-1})$,   let $a_n$ be any integer $>a_{n-1}$ and $\FF =\EE \oplus \OPN(a_n)$, so that $\gst(\FF)=(a_1,\dots,a_{n})$.  Then every  subsheaf of $\FF$ isomorphic to $\oplus_{i =1}^{n-1}\OPN (a_i)$ is strictly contained in $\EE$ and so it can not be a maximal free subsheaf of $\FF$ in the sense of \cite{MDP}.
\end{description}
\end{example}

\begin{lemma}\label{lemmafacile} 
Let $\FF$ be a rank $n$ reflexive sheaf on $\PP^N$ and let $\FF_H$ be its 
restriction to a general linear space $H \cong \PP^r$ in $\PP^N$. If $\gst 
(\FF)=(a_1, \dots ,a_n)$ and $\gst (\FF_H)=(a_1', \dots, a_n')$ (as a reflexive 
sheaf on $\PP^r$), then: 
\[a_i \leq  a_i' \hbox{ \ \  for every \ \ } i=1, \dots n.\] 
In particular, if  $char \, k=0$: \ $a_i \leq b_i$, where $(b_1, \dots, b_n)=st(\FF)=gst(\FF_L)$ 
for a general line $L$ in $\PP^N$. 
\end{lemma} 

\begin{proof}\ 
Consider an injective map $\phi :\oplus_{i=1}^n\OPN(a_i) \hookrightarrow \FF$. 
For a general linear space $H$ we have  $\mathcal{T}or^1(\coker (\phi) , 
\OO_H)=0$ \ and then the the restriction 
\[ 
\phi_L:\oplus_{i=1}^n\OO_{H}(a_i) \freccia \FF_H 
\] 
is still injective, so that  $a_i\leq a_i'$. 

Moreover, for a general line $L$, $ \FF_L\cong \oplus_i\OO_L(b_i)$, so that 
the splitting type of $\FF$ is precisely the global section type of the 
restriction $\FF_L$ to a general line $L$. 
 
\end{proof} 

\begin{lemma}\label{lemmaprinc} 
Let $\FF$ be a rank $n$ reflexive sheaf on $\PP^N$ and let 
$\gst(\FF)=(a_1,\dots,a_n)$.  If $\gsr(\FF) \geq 
c$, then: 
\begin{enumerate}[(i)] 
\item the integers $a_1, \dots, a_c$  are non-negative; 
\item every $f \colon \oplus_{i=1}^c \OPN (\alpha_i) \hookrightarrow \FF$ 
(where $\alpha_1 \geq \dots \geq \alpha_c\geq 0$) factorizes through 
\[\tilde{f} \colon \oplus_{i=1}^c \OPN (\alpha_i) \hookrightarrow \GG  \hbox{ \ 
and \ }\hat{f}:\GG\hookrightarrow \FF\] 
where $\GG$ is  a rank $c$ reflexive sheaf  and $\coker(\hat{f})$ is 
torsion free. 
\item If $char(k)=0$ and $\st(\FF)=(b_1,\dots,b_n)$,   then $b_1,\dots, b_c$ are non-negative; moreover if $b_{c+1},\dots,b_n$ are strictly negative, then 
$\mathrm{st}(\GG)=(b_1,\dots,b_c)$. 
\item If $c=n-1$, then   $\hat{f}$ can be included in the exact sequence 
\begin{equation}\label{lasequtileatutto} 
0\freccia \GG\stackrel 
{\hat{f}}{\freccia}\FF\stackrel{g}{\freccia}\IY(q)\freccia 0. 
\end{equation} 
  where $Y$ is either empty or a codimension $\geq 2$ subvariety of $\PP^N$. 
\end{enumerate} 
\end{lemma} 

\begin{proof} Part \textit{(i)} is just the note following Definition \ref{typeglobalsections}. 
\begin{enumerate} 
\item[\textit{(ii)}]  Let $\RR$ be  the rank $n-c$  sheaf $\coker (f )$. 
Dualizing  twice the exact sequence 
\[ 
0\freccia \oplus_{i=1}^c\OPN(\alpha_i)\stackrel{f}{\freccia} \FF\freccia 
\RR\freccia0. 
\] 
we obtain 
\begin{equation}\label{Rdoppioduale} 
0\freccia \GG\stackrel{\hat{f}}{\freccia} 
\FF^{\vee\vee}\stackrel{g}{\freccia} \RR^{\vee\vee} 
\end{equation} 
where $\GG$ is a rank $c$ reflexive sheaf  containing 
$\oplus_{i=1}^c\OPN(\alpha_i)$, so that 
$\gsr(\GG)=\rk(\GG)=c$ and there is $\tilde{f} \colon \oplus_{i=1}^c \OPN (\alpha_i) \hookrightarrow \GG$. Furthermore the sheaf $\coker(\hat{f})$ is torsion-free because it is a subsheaf of $\RR^{\vee\vee}$ which is reflexive. 
\item[\textit{(iii)}]  By Lemma \ref{lemmafacile} and part \textit{(i)}, we know that 
$b_1,\dots,b_c\geq 0$ and also that $\beta_1 \geq \dots \geq \beta_c \geq 0$, 
where  $(\beta_1, \dots, \beta_c)=st(\GG)$. 

Suppose $b_{c+1},\dots , b_n<0$. For a  general line $L$: 
\[ 
\GG_L\cong\oplus_{j=1}^c\OO_L(\beta_j)\stackrel{\hat{f}_L}{\hookrightarrow}\FF_L\cong\oplus_{i=1}^n\OO_L(b_i). 
\] 
So $\im(\hat{f}_L)$ is a rank $c$ subsheaf of $\oplus_{i=1}^c\OO_L(b_i)$ and then the 
quotient 
$ 
\oplus_{i=1}^{c}\OO_L(b_i)/\im(\hat{f}_L) 
$ 
is a torsion sheaf;  on the other hand, it is also isomorphic to a subsheaf of 
$(\RR^{\vee\vee})_L$ which is torsion-free. Thus 
$\GG_L\cong \im(\hat{f}_L)\cong\oplus_{i=1}^c\OO_L(b_i)$. 
\item[\textit{(iv)}] If  $n=c+1$,  the sheaf $\RR^{\vee\vee}$ given by 
(\ref{Rdoppioduale}) is a rank 1 
reflexive sheaf, that is $\RR^{\vee\vee}\cong \OPN(s)$. Then $\im 
(g)$ is a subsheaf of $\OPN(s)$; we can write it as a suitable  twist  of the 
ideal 
sheaf of a  subvariety $Y$ of codimension $\geq2$ and get  the exact sequence 
(\ref{lasequtileatutto}). 
\end{enumerate} 
 
\end{proof}

\section{Sharp lower bounds on $c_{1}$ and $c_{2}$} 

In this section $\FF$ will always be a rank $n$  reflexive sheaf on 
$\PP^N$ with Chern classes $c_i$ and global section type $\gst(\FF)=(a_1, 
\dots, a_n)$. We want to state some relations between $c_1(\FF)$ and 
the number $\delta (\FF): =\sum_{i=1}^n a_i$  and also between 
$c_2(\FF)$ and the number $\gamma(\FF):=\sum_{1\leq i<j\leq n} a_ia_j$. By 
the definition itself of type by global sections, there is a maximal injective 
map $\oplus_{i=1}^n\OPN (a_i) \hookrightarrow \FF$; the integers $\delta (\FF)$ 
and $\gamma (\FF)$ are precisely the first and second Chern class of $\oplus 
\OPN (a_i)$. When $n=1$, we assume $\gamma(\FF)=0$.

First of all, we collect some properties that we will use many times. 

\begin{lemma}\label{osservazione1}  Let $\FF$ be a rank $n$ reflexive sheaf.

\begin{enumerate} 
\item[(i)]  $c_1(\FF(l))-\delta(\FF(l))$ does not depend on $l$ that is 
$c_1(\FF(l))-\delta(\FF(l))=c_1- \delta (\FF)$; 
\item[(ii)]  $c_2(\FF(l))-\gamma(\FF(l))= c_2-\gamma (\FF)+(c_1 - 
\delta(\FF))\cdot  (n-1)\cdot  l$. 
\end{enumerate}
We now assume $a_n=0$ (and then $\gsr(\FF)=n$).
\begin{enumerate}
\item[(iii)] There is an exact sequence: 
\begin{equation}\label{dausare} 
0 \freccia \GG \freccia \FF \freccia \II_Y(q) \freccia 0 
\end{equation} 
where $\GG$ is reflexive, $\rk(\GG)=\gsr (\GG)=n-1$, $\gst(\GG)=(a_1, \dots, 
a_{n-1})$, $H^0\II_Y(q)\neq 0$, $q\geq 0$, $H^0\IY(q-1)=0$, $\gamma(\GG)=\gamma(\FF)$, $\delta(\GG)=\delta(\FF)$. 
\end{enumerate} 
Moreover: 
\begin{enumerate} 
\item[(iv)] if $q=0$, then $Y$ is empty; the converse is true with the hypothesis  $H^1\GG(-q)=0$ and in both cases $\FF\cong \GG\oplus\OPN$.
\item[(v)] If $\GG $ is a direct sum of line bundles and $q=1$, then $Y$ is a 
complete intersection $(1,r)$ where $r=c_2(\FF)-\gamma (\FF)-\delta (\FF)$ and 
$\FF$ has the short free resolution: 
\begin{equation}\label{libera} 
0\freccia \OPN(-r) \freccia \oplus_{i=1}^n \OPN (a_i)\oplus \OPN(-r+1) \freccia 
\FF \freccia 0 
\end{equation} 
\end{enumerate} 
\end{lemma} 

\begin{proof} 
Part \textit{(i)} and part \textit{(ii)} can be easily obtained by a straightforward computation. 

To show \textit{(iii)} we apply   Lemma \ref{lemmaprinc} to the map 
$\oplus_{i=1}^{n-1}\OPN (a_i) \hookrightarrow \FF$.  If $\gst(\GG)=(a_1', \dots, 
a_{n-1}')$ then  $a_i' \geq a_i$ because there is an injective map 
$\oplus_{i=1}^{n-1}\OPN (a_i) \hookrightarrow \GG$; on the other hand, 
$\oplus_{i=1}^{n-1}\OPN(a_i') \hookrightarrow \GG \hookrightarrow \FF$ implies 
$a_i' \leq a_i$. 

Moreover we  have $H^0\II_Y(q) \neq 0$ so that $q\geq 0$: otherwise 
$\GG$ would be isomorphic to the rank $n$ subsheaf of $\FF$ generated by its 
global sections. Observe that, thanks to the assumption $a_n=0$, we also have $h^0\IY(q-1)=0$.
Finally, the equalities on $\delta$ and $\gamma$ are immediate consequence of the 
assumption $a_n=0$.

\medskip 
For \textit{(iv)} observe that $q=0$ implies $Y=\emptyset$ (because $H^0\II_Y\neq 0$) so that   the sequence (\ref{dausare})
 is exact on
 global sections and splits, that is   $\FF=\GG\oplus\OPN$. On 
the other hand, if $Y=\emptyset$ and $H^1\GG(-q)=0$, then $\FF$ is the trivial extension, because $\Ext^1(\OPN(q),\GG)\cong H^1\GG(-q)$; moreover $q$ must be $0$ because 
$\gst(\FF)=(a_1, \dots, a_{n-1},q)=(a_1, \dots, a_n)$ and  $a_n=0$. 

\medskip 

Finally, in order to prove \textit{(v)}, assume $\GG\cong \oplus_{i=1}^{n-1}\OPN(a_i)$  so that $Y$ is a 
subvariety of pure codimension 2  (see Lemma \ref{spezzecodim2}); if $q=1$, then $Y$   is 
contained in a hyperplane, so that it is a complete intersection $(1,r)$, where 
$r=\deg_2(Y)$. 

Using mapping cone on (\ref{dausare}) and the standard free resolution: 
\[0\freccia \OPN(-r) \freccia \OPN (-r+1)\oplus\OPN \freccia \II_Y(1) \freccia 0\] 
we get the required   free resolution for $\FF$. 
 
\end{proof} 

\begin{remark}\label{osservazione2}  The assumption $a_n=0$ which appears in Lemma \ref{osservazione1}  will play a key role in the 
following, because very often it leads to  easier computations. For instance, if $\FF$ and 
$\GG$ are as in Lemma \ref{osservazione1}  \textit{(iii)} and $a_n=0$, then 
$\delta (\GG)=\sum_{1\leq i\leq n-1}a_i=\sum_{1\leq i\leq n}a_i=\delta (\FF)$ 
and  $\gamma (\GG)=\sum_{1\leq i<j \leq n-1}a_ia_j=\sum_{1\leq i<j \leq 
n}a_ia_j=\gamma (\FF)$: in the sequel when $a_n=0$ we will use simply $\delta$ 
and $\gamma$  and $\sum a_i$, $\sum a_ia_j$ for both sheaves. 
\end{remark} 
\begin{theorem}  \label{positconai} 
Let $\FF$ be a reflexive sheaf with $\gst(\FF)=(a_1,\dots,a_n)$. If  $\gsr (\FF)=rk(\FF)=n$, then: 
\begin{equation*} 
c_1 (\FF) \geq   \sum_{i=1}^n a_i    \hspace{1cm}, \hspace{1cm} 
c_2(\FF)\geq  \sum_{1\leq i<j\leq n}a_ia_j; 
\end{equation*} 
and equality holds in either case if and only if $\FF\cong 
\oplus_{i=1}^n\OPN(a_i)$. 
\end{theorem}

\begin{proof} Thanks to Lemma \ref{osservazione1} \textit{(i)} and \textit{(ii)}, it is sufficient to prove the statement 
for the minimal twist of $\FF$ which satisfies our hypothesis: thus   without 
lost in generality we may assume $a_n=0$. 

\medskip

We proceed by induction on $n$. As the statement clearly holds for line bundles, 
suppose  $n\geq 2$ and  the statement   true for any rank $n-1$ 
reflexive sheaf;  thus it holds in particular for  the sheaf $\GG$ defined in Lemma \ref{osservazione1} \textit{(iii)}.

For the inequality on $c_1$, we just have to observe that
\[
c_1(\FF)=c_1(\GG)+q\geq c_1(\GG)\geq \sum_{i=1}^{n-1}a_i
\]
where the last inequality is obtained applying the inductive hypothesis on $\GG$. Using the assumption $a_n=0$ (so that $\delta:=\delta(\FF)=\delta(\GG)$ and $\gamma:= \gamma(\FF)=\gamma(\GG)$), we have $c_1(\FF)\geq \delta(\FF)=\delta(\GG)$.
If $c_1(\FF)=\delta$, then also $c_1(\GG)=\delta$ so that, by the inductive hypothesis, $\GG \cong \oplus_{i=1}^{n-1}\OPN (a_i)$ and moreover $q=0$, so that $\FF\cong \GG \oplus \OPN$ (see Lemma \ref{osservazione1} \textit{(iv)}).

For  the second Chern class we have: 
\[ c_2(\FF) =c_2(\GG)+c_1(\GG) q+ \deg_2(Y) \geq \gamma \] 
where the last inequality is due to the inductive hypothesis $c_2(\GG)\geq \gamma$. 

\medskip 

Finally, if $c_2(\FF)=\gamma$,   we also have   $c_2(\GG)=\gamma $ (so that by induction 
$\GG\cong \oplus_{i=1}^{n-1}\OPN(a_i)$) and $\deg_2(Y)=0$ (so that again 
$Y=\emptyset$, thanks to Lemma \ref{spezzecodim2}). Observing that in this case $H^1_*\GG=0$, we can conclude using lemma \ref{osservazione1} \textit{(iv)}. 
 
\end{proof} 

\begin{remark}
The result of Theorem \ref{positconai}  on $c_1(\FF)$ holds without the hypothesis $\gsr (\FF)=rk(\FF)=n$, because $c_1(\FF)-\sum a_i$ is invariant by twist, while this hypothesis is necessary for the result about $c_2(\FF)$ because $c_2-\sum a_ia_j$ is not invariant by twist (Lemma \ref {osservazione1}, $(ii)$).
\end{remark}

The following example shows that for every choice of non-negative integers 
$a_1, \dots, a_n$ and $s$, there are rank $n$ reflexive sheaves $\FF$ such that 
$\gst(\FF)=(a_1, \dots, a_n)$ and $c_1(\FF)=\sum a_i +s$; in other words, $c_1(\FF)$ can in fact reach every value  above the minimal one $\delta(\FF)$  given by the previous 
results. 

\begin{example}\label{esempio1} Let $a_1\geq  \dots \geq  a_n$  and $s$ be non-negative integers. We define $p\geq a_1-a_2+ s$ and let $Y$ be a 
complete intersection $(s,p)$ in $\PP^N$. 

As $\omega_Y\cong \OO_Y(s+p-N-1)$, then  $\omega_Y(N+1 -s 
-a_1+a_2)\cong\OO_Y(p-a_1+a_2)$ has a section which generates it almost 
everywhere. Such a section gives  an extension: 
\[ 0 \rightarrow \OPN (a_2) \rightarrow \GG \rightarrow \II_Y(s +a_1) 
\rightarrow 0 \] 
where $\GG$ is a rank $2$ reflexive sheaf with first Chern class $c_1 
(\GG)=a_1+a_2+s $; moreover  by the hypothesis on $p$ we have $H^0\II_Y(s-1)=0$ and $H^0\II_Y(s+a_1-a_2-1)=H^0\II_Y(s)\otimes H^0\OPN(a_1-a_2-1)$, so that  $\gst (\GG)=(a_1,a_2)$. 

Finally, the rank $n$ reflexive sheaf  $\FF=\GG \oplus \OPN (a_3) \oplus \dots 
\oplus \OPN (a_n)$ has first Chern class $c_1 (\FF)=\sum a_i+s$ and  
global section type $\gst (\FF)=(a_1,\dots, a_n)$. 

Note that such a sheaf has the short free resolution: 
\begin{equation}\label{disegusuc2} 
0\freccia \OO_{\PP^N}(a-s)\freccia  \oplus_{i=1}^n\OO_{\PP^N}(a_i)\oplus\OO_ 
{\PP^N}(a)\freccia \FF\freccia 0 
\end{equation} 
with $a=a_1-p+s$.
\end{example} 

All the  sheaves obtained in Example \ref{esempio1} are of a very special type: 
they have homological dimension $\leq 1$ and a  short free resolution. Now we 
will show that in fact every sheaf for which $c_1(\FF)=\delta(\FF) +1$ has the 
same nice properties. Note that we can obtain reflexive sheaves $\FF$ with a resolution of the type (\ref{disegusuc2}) and $\gst (\FF)=(a_1,\dots, a_n)$  for every possible  integer   $a \leq  a_2 $. 

\begin{corollary}\label{c1+1} Let $\FF$ be a reflexive sheaf.  If $\FF$ is not a 
sum of line bundles, then \[c_1(\FF)\geq \sum a_i +1\] and equality holds if and only 
if   $\hd(\FF)\leq 1$  and $\FF$ has a free resolution (\ref{libera}). 
\end{corollary} 

\begin{proof} We   have to show that $c_1(\FF)-\delta(\FF) \geq 1$. As 
$c_1(\FF)-\delta (\FF)$ does not change up a twist, we may assume $a_n=0$: so the inequality follows from Theorem \ref{positconai}. 

\medskip 

Moreover $a_n=0$ also allows us to use the exact sequence (\ref{dausare}). 
Assume that $c_1(\FF)=\delta(\FF) + 1$. If $q=0$, then $\FF\cong \GG\oplus \OPN$ (Lemma \ref{osservazione1}, \textit{(iv)}) and we 
conclude by induction on the rank.  If $q\neq 0$, then  $q=1$ and   we have 
$c_1(\GG)=c_1(\FF)-q= \delta(\FF)=\delta(\GG)$, so that  $\GG \cong 
\oplus_{i=1}^{n-1} \OPN (a_i)$ (see Theorem \ref{positconai}) and we conclude thanks to Lemma \ref{osservazione1} \textit{(v)}. 
 
\end{proof} 

\medskip

For what concerns the second Chern class,  not all values    above the minimal 
one given in Theorem \ref{positconai} can   in fact be realized, at least for 
general $a_1, \dots, a_n$.  For instance, if $a_2>0$, any integer   in the 
interval $  [ \gamma(\FF)+1 \ , \ \gamma (\FF)+\delta (\FF) -a_2  ]$        is
not allowed. More precisely:

\begin{theorem}\label{miglioramento} 
Let $\FF$ be  rank $n\geq 2$ reflexive sheaf on $\PP^N$ with Chern classes $c_i$ and global section type $(a_1, \dots, a_n)$.  Then:
\begin{equation}\label{piuforte} c_2\geq \sum_{1\leq i<j\leq n}a_ia_j+\left( c_1-\sum_{1\leq i\leq n}a_i\right) \left( 
\sum_{1\leq i\leq n}a_i+1 -a_2\right).\end{equation} 
\end{theorem} 

\begin{proof} First of all observe that  equality holds in (\ref{piuforte}) if $\FF \cong \oplus \OPN(a_i)$ is a sum of line bundles, because in this case  $c_1(\FF)=\sum a_i$ and  $c_2(\FF)=\sum a_ia_j$.  So, assume that  $\FF$ is not split.

\medskip

We have to show  that the integer $\Delta(\FF)=c_2-  \gamma - ( c_1-\delta )  ( 
\delta+1 -a_2 )$ is non-negative.  Using Lemma \ref{osservazione1} \textit{(i)} and \textit{(ii)}, we can see that $\Delta(\FF)$ is invariant under twist, that is  $\Delta(\FF(l))=\Delta(\FF)$ for every $l\in \ZZ$.  
Thus,  without lost in generality,  we may assume $a_n=0$.

\medskip 

Let $\GG$ and $Y$ be as   in    Lemma \ref{osservazione1} \textit{(iii)} and denote by $d_i$ the Chern classes of $\GG$. 
By the exact sequence (\ref{dausare}) and Lemma \ref{spezzecodim2}, we get:
\[c_2 = d_2+(c_1-d_1)d_1+\deg_2(Y) \geq \gamma + ( d_1-\delta )  ( 
\delta+1 -a_2 ) +(c_1-d_1)d_1+\deg_2(Y). \] 
If $\rk(\FF)=2$, then $\GG\cong \OPN(a_1)$ so that the exact sequence (\ref{dausare}) becomes
\[
0 \freccia \OPN(a_1)\freccia \FF \freccia \II_Y(c_1-a_1) \freccia 0 
\]
 and $Y$ is a subvariety of pure codimension $2$ (see Lemma \ref{spezzecodim2}) not contained in hypersurfaces of degree  $c_1-a_1-1$, because on the contrary  $a_2\geq 1$ against the assumption. Moreover  $d_1=\delta=a_1$,  $d_2=\gamma=0$. 
So it will be sufficient  to prove  that $(c_1-d_1)d_1+\deg_2(Y)\geq (c_1-a_1)(a_1+1)$ that is $\deg_2(Y)\geq (c_1-d_1)$. This is true because every subvariety $Y$ of pure codimension $2$ and degree $s$ is always contained in  some degree $s$ hypersurfaces (for instance  cones). 

\medskip

If $\rk(\FF)\geq 3$, we can proceed by induction on the rank and assume that (\ref{piuforte}) holds for the reflexive sheaf $\GG$. Thanks to the equality:
$$\Delta(\FF)=\Delta(\GG)+\deg_2(Y)-(c_1-d_1)(\delta +1-a_2-d_1)$$
and the inductive hypothesis,  it is sufficient to prove that:
\begin{equation}\label{boundc2}
\deg_2(Y)\geq (c_1-d_1)(\delta +1-a_2-d_1).
\end{equation}
We know that $c_1\geq d_1\geq \delta$ and $a_2\geq a_n=0$ (see Theorem \ref{positconai}); then (\ref{boundc2}) clearly holds if  either $a_2>0$ or $d_1>\delta$.

The only case left to consider is $a_2=d_1-\delta=0$; in this case  $\GG\cong  \OPN(a_1) \oplus \OPN^{n-2}$ (again by Theorem \ref{positconai}) so that $d_1=a_1$,  $Y$ is a subvariety of pure codimension $2$ (see Lemma \ref{spezzecodim2}),  $H^0\II _Y(c_1-a_1-1)=0$ and  the exact sequence (\ref{dausare}) becomes
\[
0 \freccia \OPN(a_1)\oplus \OPN^{n-2}\freccia \FF \freccia \II_Y(c_1-a_1) \freccia 0. 
\]Thus  we can now conclude  as in the rank $2$ case. 
 
\end{proof} 

  Now we list a few remarkable consequences of Corollary \ref{c1+1} and Theorem \ref{miglioramento}. Note that  the integer $\sum a_i-a_2+1=\sum_{i\neq 2}a_i +1$ is strictly positive when  $\gsr(\FF) =n$.
 
\begin{corollary} \label{boundperc2}  Let $\FF$ be a rank $n$ reflexive sheaf and 
 let $\alpha_1,  \dots , \alpha_n $ be integers such that $\FF$ has a free subsheaf isomorphic to $ \oplus_{i=1}^n\OPN (\alpha_i) $: for instance $(\alpha , \dots , \alpha_n)=\gst(\FF)$. 
  If  $\alpha_1 \geq \dots \geq \alpha_n \geq 0$, then:
\begin{equation*}
	c_1\geq   \sum_{1\leq i\leq n} \alpha_i 
	\hspace{.5cm},\hspace{.5cm}
	c_2\geq   \sum_{1\leq i<j\leq n}\alpha_i\alpha_j.
\end{equation*}
Moreover equality on $c_1 $ holds  if and only if   $\FF$ is  $\oplus_{i=1}^n
\OPN (\alpha_i)$; equality on $c_2$ holds if and only if $\FF$ is either
$\oplus_{i=1}^n \OPN (\alpha_i)$ or $\OPN ^{n-1}\oplus \OPN (a_1)$ for some $a_1
>\alpha_1 $.

\medskip

If $\FF$ is not a direct sum of line bundles then
\begin{equation} \label{senzaa1}
c_2(\FF)\geq \sum_{1\leq i<j\leq n}\alpha_i\alpha_j +\sum_{i=1}^n \alpha_i+1-\alpha_2
\end{equation}
and moreover
\begin{equation} \label{senzaa2}
c_2(\FF)\geq \sum_{1\leq i<j\leq n}\alpha_i\alpha_j +2(\sum_{i=1}^n \alpha_i+1-\alpha_2)
\end{equation}
unless $\FF$ has the short free resolution given in (\ref{libera}).\end{corollary}

\begin{corollary}\label{positivi}
Let $\FF$ be a reflexive sheaf such that  $\gsr(\FF)=\rk(\FF)$. Then
\[
c_1(\FF)\geq 0 \qquad \text{,}\qquad c_2(\FF)\geq 0
\]
and either equality holds only if $\FF =\OPN^{n-1} \oplus\OPN(c_1)$. 
\end{corollary}

\begin{example} Sharp cases for the lower bounds on $c_2$ given by  the previous results when $(\alpha_1, \dots, \alpha_n)=(a_1, \dots, a_n)=\gst(\FF)$ can be found  in 
Example \ref{esempio1} for special values of the  parameters. 

For every $a_1\geq \dots\geq a_n\geq 0$, we can choose  $s=1$ and $p=a_1-a_2+1$ and get a sheaf $\FF$ such that $c_1=\delta +1$ and $c_2=\gamma +(\delta -a_2+1)$.

For what concerns  sheaves with homological dimension $\geq 2$ (and then without the free resolution (\ref{libera})) a sharp case is given by the rank $2$ reflexive sheaf  $\FF$ that we will consider in Example \ref{2piani}: its  global section type is $(a_1=0, a_2=0)$  and  its first and second Chern classes  are
$c_1=2$ and $c_2=2$.  Using such a  sheaf $\FF$ we can also find  reflexive sheaves of any rank $n> 2$ on $\PP^4$ realizing the equality in (\ref{senzaa2}):  for instance $\FF(l)\oplus\OO_{\PP^4}(a_3)\oplus\cdots\oplus\OO_{\PP^4}(a_n)$, where $l\geq a_3\geq \dots\geq a_n\geq 0$.
\end{example}

We conclude this section  showing that in the previous results we cannot simply avoid  the assumption on the global section rank in order to get   the positivity  of $c_1$ and $c_2$.

\begin{example}
Let $q$  be an integer, $q>>0$, and $\EE$ be the rank 2 bundle on $\PP^3$ defined as an extension
\[
0\freccia \OP\freccia \EE\freccia \IY(-2q+4)\freccia 0
\]
where $Y$ is a $(-2q)$-subcanonical double structure on a line obtained by Ferrand construction (see \cite{HSVB}, Theorem 1.5).

We have 
\[
c_1(\EE)=-2q+4 \qquad c_2(\EE)=2 \qquad  h^0\EE\neq 0.
\]
Then for every integer $t$, $0<t<q-2$,  the global section rank is  strictly lower than the rank (in fact  $\gsr(\EE(t))=\rk(\EE(t))-1$)  and both 
$
c_1(\EE(t))$ and $ c_2(\EE(t))$ are strictly negative.
\end{example}

\section{Special splitting types} 

As we have  seen in the previous section,  not all sheaves whose global section rank is lower than the rank have positive $c_1$ and $c_2$, not even vector bundles. In this section we assume something weaker about the global section rank (namely $\gsr(\FF)\geq \rk(\FF)-1$), while we introduce balancing hypothesis  on the splitting type.

\begin{proposition}\label{unasezinmeno} 
Let $\FF$ be a reflexive sheaf on $\PP^N$ such that $\rk(\FF)=n$, $\gsr(\FF)=n-1$, $c_1\leq 0$ and 
$\st(\FF)=(0,\dots,0,c_1)$.   
Then: \[c_2(\FF)\geq 0.\] Moreover  $c_2(\FF)=0$ if and only if $\FF\cong \OPN(c_1)\oplus 
\OPN^{n-1}$. 
\end{proposition} 

\begin{proof} 
Let $\GG$ and $Y$ be as in Lemma \ref{lemmaprinc}: under  the present  assumption on $\FF$,  we have $\rk(\GG)=\gsr(\GG)=n-1$ and    $\st(\GG)=(0,\dots,0)$ so that $c_1(\GG)=0$. The only sheaf of such a type is $\GG \cong \OPN^{n-1}$ (see Corollary \ref{positivi})  and so  by (\ref{lasequtileatutto}) we get $c_2(\FF)=c_2(\II_Y(c_1))=\deg_2(Y)\geq 0$. 

Moreover  $Y$ has pure codimension 2 or it is empty (see Lemma \ref{spezzecodim2}); then
 $c_2=0$ if and only if $Y=\emptyset$ and  $\FF\cong \OPN(c_1)\oplus \OPN^{n-1}$. 
 
\end{proof} 

\begin{remark} Most properties  concerning the   sheaf $\FF$ that appear in the previous section are in general \lq\lq stable by positive twist\rq\rq ; for instance for every  $l\geq 0$, if $c_1(\FF) \geq 0$ and $c_2(\FF) \geq 0$, then also $c_1(\FF(l)) \geq 0$ and $c_2(\FF(l)) \geq 0$ and if $\gsr(\FF)=\rk(\FF)$ then also $\gsr(\FF(l))=\rk(\FF(l))$. On the contrary neither the  hypothesis nor the thesis that appear in Proposition \ref{unasezinmeno} are  \lq\lq stable\rq\rq, as shown by the following example.
\end{remark}

\begin{example} 
Let $Y$ be a line in $\PP^3$. For any integer 
$c\leq  -2$, $\omega_Y(4-c)\cong\OO_Y(2-c)$ has a section which generates it 
almost everywhere.  Therefore there is an extension (see \cite{HSRS}, proof 
of 1.1) 
\[ 
0\freccia \OO_{\PP^3}\freccia \EE\freccia \II_Y(c)\freccia 0 
\] 
which defines   the   rank 2 reflexive sheaf $\EE$.  
This  sheaf  satisfies the hypothesis of Proposition 
\ref{unasezinmeno} that is $\gsr(\EE)=1=\rk(\EE)-1$,  $c_1(\EE)=c<0$, $\st(\EE)=(0,c_1)$;  as a consequence, the second Chern class 
$c_2(\EE)$ is positive (in fact $c_2(\EE)=\deg_2Y=1$). 
On the other hand, $\EE(1)$ has a negative first Chern class, but it does not satisfy the hypothesis of  Proposition \ref{unasezinmeno}  about the splitting type, which is in fact   $(1, c+1)$ instead of $(0, c+2)$. If we compute its second Chern class $c_2(\EE(1))=c+2$, we can see that $c_2(\EE(1))$ is strictly negative  when $c<-2$ and  moreover that  $c_2(\EE(1))=0$ 
    when $c=-2$ although $\EE(1)$ is not a split bundle. 
\end{example}

\begin{proposition} 
Let $\FF$ be a reflexive sheaf on $\PP^N$ with Chern classes $c_i$, such that $\rk(\FF)=n$, $\gsr(\FF)=n-1$, $c_1\leq 0$ and 
$\st(\FF)=(1,0,\dots,0,c_1-1)$. Then $c_2\geq c_1-1$. Moreover\\ 
\begin{enumerate}
\item[(i)] $c_2=c_1-1$ if and only if $\FF\cong \OPN(1)\oplus\OPN(c_1-1)\oplus\OPN^{n-2}$.
\item[(ii)] If $H^0\FF(-1)\neq 0$, then $c_2=c_1$ if and only if $\FF$ has the following free resolution:
\begin{equation}\label{laurea}
0\freccia \OPN(c_1-3)\freccia \OPN(c_1-2)^2\oplus\OPN(1)\oplus\OPN^{n-2}\freccia \FF\freccia 0.
\end{equation}
\item[(iii)]If $H^0\FF(-1)=0$ and furthermore $H^1\FF_{H_r}(-1)=0$ for every general linear subspace $H_r\cong \PP^r$  in 
$\PP^N$ ($r\geq 3$), then $c_2=c_1$ if and 
only if $\FF\cong \GG\oplus \OPN(c_1-1)$ and $\GG$ has the free resolution:
\begin{equation}\label{liberaperG} 
0\freccia \OPN(-1) \freccia \OPN^{n} \	 \freccia 
\GG \freccia 0 
\end{equation} 
\end{enumerate}
\end{proposition}

\begin{proof} 
Let $\GG$ be the reflexive sheaf constructed in Lemma \ref{lemmaprinc} (\ref{lasequtileatutto});  denote by $d_i$ its Chern classes.  We have $\rk(\GG)=\gsr(\GG)=n-1$ (so that $d_2 \geq 0$: see Corollary \ref{positivi}) and moreover 
$\st(\GG)=(1,0,\dots,0)$ (so that $d_1=1$).  Using the exact sequence (\ref{lasequtileatutto}) we obtain:
\begin{equation}\label{disugc2-c1+1}
c_2=d_2+d_1(c_1-d_1)+\deg_2(Y)=d_2+c_1-1+deg_2(Y) \geq c_1-1.
\end{equation}
Note that $c_1-d_1 =c_1-1<0$; then $H^0\II_Y(c_1-1)=0$ so that $H^0\FF(-1) \cong H^0\GG(-1)$.

\medskip

\textit{(i)}
The equality $c_2=c_1-1$ can hold only if $d_2=0$ so that $\GG \cong \OPN (1)\oplus \OPN^{n-2}$ (Corollary \ref{positivi}); moreover $Y$ is a pure codimension $2$ subvariety of degree $0$ that is $Y=\emptyset$ and $\FF\cong \OPN (1)\oplus \OPN^{n-2}\oplus \OPN(c_1-1)$.

\medskip

The next value $c_2=c_1$  can be realized if and only if either $d_2=0$ and $\deg_2(Y)=1$ or $d_2=1$ and $\deg_2(Y)=0$. 

\medskip

 \textit{(ii)} In the first case $\GG \cong \OPN (1)\oplus \OPN^{n-2}$, $Y$ is a pure codimension $2$ subvariety of degree $1$, that is a complete intersection $(1,1)$.
   Using  mapping cone on (\ref{lasequtileatutto}) and the canonical free resolution of $\IY(c_1-1)$, we  get (\ref{laurea}).

\medskip

\textit{(iii)} In the other case, namely  if $d_2=1$ and $\deg_2(Y)=0$,   we find  $\gst(\GG)=(0, \dots, 0)$ and then $\GG$ has the free resolution (\ref{liberaperG}) (see Corollary \ref{boundperc2}) and so $H^1_*\GG=0$.  To complete the proof we just have to show that $Y$ 
is empty. 

If not, let $r$ be  the codimension of $Y$, $3\leq r \leq N$;  the restriction  $Y'=Y\cap H_r$ to a sufficiently general linear subspace $H_r$  is a finite set of points, whose degree $\delta=\deg(Y)$ is given by  
$h^0\OO_{Y'}(c_1-2)=h^1\II_{Y'}(c_1-2)$. On the other hand $h^1\II_{Y'}(c_1-2)\leq h^1\FF_{H_r}(-1)+h^2\GG(-1)=0$. Then $Y$ is 
empty and $\FF \cong \GG \oplus \OPN(c_1-1)$ (see Lemma \ref{osservazione1}, \textit{(iv)}).
 
\end{proof}

\section{Positivity and border cases for $c_3$} 

In  the present, last section we want to investigate general conditions that ensure the positivity of the third Chern class $c_3$ of a reflexive sheaf $\FF$ of any rank. 

To this aim it will be useful to have  a deeper knowledge  about the singular locus of  $\FF$ or, more generally, about the set of points where it is not freely generated.

\begin{proposition} \label{codimluogosing}
Let $\FF$ be a rank $n$ reflexive sheaf on $\PP^N$ such that 
$\hd(\FF)\leq 1$,  $\hd(\FF^{\vee})\leq 1$.

 Then the singular locus $\mathcal S$ of  $\FF$ 
is either a  codimension  $3$ closed subset or it is empty (that is $\FF$ is a vector bundle). 
\end{proposition} 

\begin{proof} Of course there is nothing to prove if $N\leq 3$. Moreover it is equivalent to prove the statement for $\FF$ or for its dual  $\FF^{\vee}$ or for some of their twists   $\FF(k)$ or $\FF^{\vee}(k) $ ($k \in 
\ZZ$). 

Assume $\mathrm{codim} (\mathcal{S})\geq 4$ and consider  a locally free resolution of $\FF$: 

\begin{equation}\label{prima} 
0\freccia \EE_1 \freccia \EE_0\freccia \FF \freccia 0 
\end{equation} 
We may assume that $\EE_0$ is a direct sum of line bundles. 
Dualizing, we get: 
\[ 0 \freccia \FF^{\vee} \freccia \EE_0^{\vee} \freccia \EE_1^{\vee}\freccia 
\mathcal{E}xt^1(\FF,\OO_{\PP^N})\freccia 0\] 
which splits in two short exact sequences: 
\begin{eqnarray}\label{seconda} &&0 \freccia \FF^{\vee} \freccia \EE_0^{\vee} 
\freccia \LL\freccia 0 \\  \label{terza} && 0 \freccia \LL   \freccia 
\EE_1^{\vee}\freccia \mathcal{E}xt^1(\FF,\OO_{\PP^N})\freccia 0 
\end{eqnarray} 
Observe that when the singular locus $\mathcal{S}$ is a finite set of points, then 
$h^i(\mathcal{E}xt^1(\FF,\OO_{\PP^N}))=0$ for every $i\geq 1$ because the 
support of $\mathcal{E}xt^1(\FF,\OO_{\PP^N})$ is contained in $\mathcal{S}$. In this case 
using  (\ref{prima}),(\ref{seconda}),(\ref{terza}) and duality on $\EE_1$ we 
find: 
\begin{equation} 
\label{quarta} 
H^0\EE_1^{\vee} \freccia H^0\mathcal{E}xt^1(\FF,\OO_{\PP^N}) \freccia H^1\LL 
\freccia H^1\EE_1^{\vee} \freccia 0 
\end{equation} 
and moreover 
\[ h^1\LL= h^2\FF^{\vee} \hbox{ and } 
h^1\EE_1^{\vee}=h^{N-1}\EE_1(-N-1)=h^{N-2}\FF(-N-1). \] 
Then $h^2\FF^{\vee} \geq h^{N-2}\FF(-N-1)$. 

We split the proof in three steps.
\noindent \begin{itemize} 
\item  First  we prove the statement in $\PP^4$. 
By hypothesis, $\mathcal{S}$ is a $0$-dimensional set and so $h^2\FF^{\vee} \geq 
h^2\FF(-5)$. 
Since the same inequality holds for every twist of $\FF$ and $\FF^{\vee}$, then  equality holds, so that $h^1\LL=h^1\EE_1^{\vee}$.  On the other hand, up a suitable twist of $\FF$, we 
have  $H^0\EE_1^{\vee}=0$ and then $H^0\mathcal{E}xt^1(\FF,\OO_{\PP^4})=0$, so that $\mathcal{E}xt^1(\FF,\OO_{\PP^N})=0$ (because it is a constant sheaf  supported on a finite set of points) and $\LL=\EE_1^\vee$. Thus we have: 
\[0 \freccia \FF^{\vee} \freccia \EE_0^{\vee} \freccia \EE_1^{\vee} 
\freccia 0\]
 and $\FF^{\vee}$  is locally free. 

\item Now we prove the statement for every $N\geq 5$ assuming that $\mathcal{S}$ is a finite set of points. As 
above,  $\mathcal{E}xt^1(\FF,\OO_{\PP^N})$ has a finite support  and, using 
again (\ref{prima}),(\ref{seconda}),(\ref{terza}) and duality on $\EE_1$ we get: 
\begin{eqnarray*} 
h^2\FF=h^3\EE_1=h^{N-3}(\EE_1^{\vee}(-N-1))= \hbox{ (note that  } N-3\geq 2 \, ) 
&&    \\ = h^{N-3}\LL(-N-1)& =h^{N-2}\FF^{\vee}(-N-1) & 
\end{eqnarray*} 
Since the same equality holds for every twist of $\FF$ and $\FF^{\vee}$, then $H_*^2(\FF)$ and $H_*^2(\FF^{\vee})$ are finite modules.   Then for every $t\gg0$, we have $h^0\EE_1^{\vee}(-t)=0$ and $h^1\LL(-t)=h^2\FF^{\vee}(-t)=0$, so that by the cohomology exact sequence of (\ref{terza}) we obtain that $H^0\mathcal{E}xt^1(\FF,\OO_{\PP^N})(-t)=0$. This implies $\mathcal{E}xt^1(\FF,\OO_{\PP^N})=0$, because it is a constant sheaf  supported on a finite set of points; from this, as above, we deduce that   $\FF$ is locally free.

\item Finally we  consider the general case and proceed by induction on $N$. Let 
$N\geq 5$ and assume the statement true for $N-1$. If $H$ is a general hyperplane,  we can apply Lemma \ref{generale} and see that  the restriction $\FF_H$  is a  reflexive sheaf on $H$, whose singular locus   $\mathcal{S}(\FF_H)=\mathcal{S}\cap H$ has codimension $\geq 4$ in $H$; moreover short locally free resolutions of  $\FF$ and $\FF^{\vee}$ restrict to short locally free resolution of  $\FF_H$ and $(\FF^{\vee})_H\cong (\FF_H)^{\vee}$ ((Lemma \ref{generale} \textit{(iii)}). Thus $\FF_H$ satisfies the same condition as $\FF$. By the inductive hypothesis,  $\FF_H$ is locally free on $H$, namely $\mathcal{S}(\FF_H)=\mathcal{S}\cap H=\emptyset$ (Lemma \ref{generale} \textit{(iv)}).   This can happen only if  $\mathcal{S}$ is at most  a finite set of points and we conclude thanks to 
the previous item. 
\end{itemize} 
 
\end{proof}

\begin{theorem}\label{positivc3} 
Let $\FF$ be a rank $n$ reflexive sheaf on $\PP^N$ generated by global sections outside a closed subset of codimension $\geq 3$. Then: 
 \[c_3(\FF)\geq 0\]
and equality $c_3(\FF)=0 $ can hold only if  $N\geq 3$ and  $\FF_H$ is a vector bundle
for every   general linear subspace $H\cong \PP^3$   in $\PP^N$.

 If moreover $\hd(\FF)\leq 1$ and $\hd(\FF^\vee)\leq 1$, then $
c_3(\FF)=0 $ only if $\FF $ is a vector bundle having  a direct summand $ \OPN^r$, for some $r\geq n-2-h^1\FF(-c_1)$.
\end{theorem} 

\begin{proof} 
 Assume that $\FF$ is not a free bundle.
Then  $n-1$ general sections of $\FF$  degenerate on a 
generically smooth codimension 2 subvariety $Y$ given by the exact sequence:
\begin{equation}\label{seqesokonek} 
0\freccia \OO_{\PP^N}^{n-1}\freccia \FF\freccia \IY(c_1)\freccia  0 
\end{equation} 
 (see Lemma \ref{smooth} and \cite{okted}, \S 2); thanks to Lemma \ref{spezzecodim2}, we can see that $Y$  has no embedded or isolated components of codimension $\geq 3$.
 
If we apply the functor $\mathcal{H}om(\cdot,\OPN)$ to the exact sequence 
(\ref{seqesokonek}), we find:
\begin{equation}\label{okonekduale}
0\freccia \OPN(-c_1)\freccia  \FF^\vee\freccia \OPN^{n-1}\freccia 
\mathcal{E}xt^1(\IY(c_1),\OPN)\freccia 
\mathcal{E}xt^1(\FF,\OPN)\freccia 0 
\end{equation} 
where $\mathcal{E}xt^1(\IY(c_1),\OPN)\cong \omega_Y(N+1-c_1)$ (see 
\cite{alggeom} Ch. III, Proposition 7.5) and  $\mathcal{E}xt^1(\FF,\OPN)$ is  supported  on $\mathcal S$.

Then $\OPN^{n-1}\freccia \omega_Y(N+1-c_1)$ is surjective outside $\mathcal S$ 
that is $\omega_Y(N+1-c_1)$ is generated almost everywhere by its global sections; as $\omega_Y(N+1-c_1)$ has rank 1, it is in fact generated almost everywhere by just one  of them. 
Such a section and the isomorphisms:
 \begin{multline*}
 H^0\omega_Y(N+1-c_1)\cong \Hom_Y(\OO_Y,\omega_Y(N+1-c_1))\cong\\
 \cong\Ext ^2_{\PP^N}(\OO_Y,\OO_{\PP^N}(-c_1))\cong \Ext ^1_{\PP^N}(\II_Y,\OO_{\PP^N}(-c_1))
\end{multline*}
(see \cite{alggeom} III, Lemma 7.4), give an extension: 
\begin{equation}\label{seqeshart} 
0\freccia \OPN\freccia \mathcal E\freccia \IY(c_1)\freccia 0 
\end{equation} 
 where $\mathcal E$  is a rank 2 reflexive sheaf.
Using multiplicativity of   Chern classes in
  (\ref{seqesokonek}) and (\ref{seqeshart}), we get
\[ 
c_3(\FF)=c_3(\IY(c_1))=c_3(\mathcal E)=c_3(\EE_H)\geq 0. 
\] 
 for every  general $H\cong \PP^3$ in $\PP^N$ (see Lemma \ref{generale} \textit{(v)} and \cite{HSRS}, Theorem 4.1). Note that our hypothesis on $\FF$ also holds for every restriction of $\FF$ to a general linear subspace.

Finally, $c_3(\FF)=0$ if and only if $c_3(\EE_H)=0$ , that is $\EE_H$ is a vector bundle  (see again   \cite{HSRS}, Theorem 4.1).  
The curve $C=Y\cap H$ is $(c_1-4)$-subcanonical that is $\omega_C(4-c_1)\cong \OO_C$ and a section which generates it almost everywhere in fact generates it;  in the exact sequence  (\ref{okonekduale})  the map $\OO_{\PP^3}^{n-1}\freccia  \mathcal{E}xt^1(\II_C(c_1),\OO_{\PP^3})\cong \OO_C$ is surjective so that $\mathcal{E}xt^1(\FF_H,\OO_{\PP^3})=0$ that is $\FF_H$ is locally free. 

    If moreover $\hd(\FF)\leq 1$ and $\hd(\FF^\vee)\leq 1$, by Proposition \ref{codimluogosing} we can conclude that $\FF$ is a vector bundle too.
   
  The exact sequence  (\ref{okonekduale}) becomes
   \[
0\freccia \OPN(-c_1)\freccia  \FF^\vee\freccia \OPN^{n-1}\freccia 
\OO_Y \freccia 0. 
\]
By a suitable change of base we can assume that $r\geq (n-2-h^0\OO_Y)$ copies of $\OPN$ are in fact contained in the kernel of the last map and $\FF$ contains a direct summand $\OPN^{r}$.
 
\end{proof}

\begin{corollary}\label{semprec3}
Let $\FF$ be a rank $n$ reflexive sheaf on $\PP^N$ ($N\geq 4$) such that $\hd(\FF)\leq 1$ and $\hd(\FF^\vee)\leq 1$. If $\FF$ is generated by global sections, then:
\[
c_3(\FF)=0 \Longleftrightarrow \FF \cong \GG \oplus \OPN^{n-2}, \hbox{ where $\GG$ is a rank $2$ vector bundle}.\]
\end{corollary}
\begin{proof} If $\FF \cong \GG \oplus \OPN^{n-2}$ and  $\GG$ is a rank $2$ vector bundle, then $c_3(\FF)=c_3(\GG)=0$. The converse is a straightforward consequence of the previous result. In fact for $n-1$ general global  sections of $\FF$, the zero locus $Y$ is a codimension 2 smooth  subvariety in $\PP^N$; as  $N\geq 4$,  $Y$ must be irreducible so that $h^1\FF(-c_1)=h^1\II_Y=h^0\OO_Y=1$.
 
\end{proof}

\begin{corollary}\label{c3nullafibrato} 
Let $\FF$ be a rank $2$ reflexive sheaf on $\PP^N$ with $\hd(\FF) \leq 1$. Then: 
\[   c_3=0  \Longleftrightarrow \FF \hbox{ is a vector bundle.} \] 
\end{corollary}

\begin{proof} For rank 2 reflexive sheaves, duality $\FF^{\vee} \cong \FF(-c_1)$ implies  $\hd(\FF)=\hd(\FF^{\vee})$; moreover the third Chern class is invariant under twist. Thus we can apply the previous result to a suitable  twist of $\FF$  and conclude.
 
\end{proof}

Note that in the previous results we cannot simply avoid either  hypothesis on the homological dimension or on global sections of $\FF$, because in fact not every  reflexive sheaf whose third Chern class vanish is a vector bundle. 

\begin{example} Let  $\FF$ be   a rank 2 reflexive sheaf with $c_2=r$ and $c_3=r^2$ (for instance we can find such a sheaf for special values   $n=2$ and $a_1=a_2=0$ in (\ref{libera}) ). Then $\FF \oplus \OPN (-r)$ is not a bundle even if $c_3=0$.
\end{example}

\begin{example}\label{2piani} 
Let $Y$ be the union of two general 2-planes in $\PP^4$ with only one common 
point $Q$. A general    non-zero global section of the sheaf $\omega_Y(3)$ 
generates it at every point  $P\in Y$ except at the point $Q$. Using    such a 
section we can define an extension: 
\[0 \freccia \OO_{\PP^4} \freccia \FF \freccia \II_Y(2) \freccia 0 \] 
where $\FF$ is a rank 2 reflexive sheaf   with $c_1=2$ and  $c_2=2$. If we cut 
by a general hyperplane $H$, we find  that  $Y_H$ is the disjoint union of two 
lines in $H=\PP^3$ and that   $\FF_H(-1)$ is a null correlation bundle. Then 
$c_3(\FF)=c_3(\FF_H)=0$  even if $\FF$ is not locally free. 

Note that in fact $Y$ is not locally Cohen-Macaulay (at the point $Q$), while 
its general hyperplane section is. If $k \gg 0$, a general global section of 
$\FF(k)$ degenerates on an integral  surface which is not locally 
Cohen-Macaulay, while its general section is a subcanonical curve in $\PP^3$. 
\end{example}

We conclude by showing that we can not expect to control the 
positivity of the $i$-th Chern class for every    reflexive sheaf when  
$i\geq 4$, even if $i$ is lower than the rank. 
In fact  for every pair $(n,i)$ 
(both $\geq 4$)  there are    rank $n$ reflexive sheaves on $\PP^N$, generated 
by global sections,  having  negative $c_i $. In the following example we construct sheaves of this type; it is not difficult to generalize the construction in order to obtain for every $t\in \mathbb{N}$ reflexive sheaves $\FF$ with the same properties and such that moreover $\FF(-t)$ is generated by global sections.
  
\begin{example}\label{negativo}
Let $\GG'$ be any rank 2 reflexive sheaf   on $\PP^N$ with third Chern class 
strictly positive.  If   $l\gg 0$, the 
sheaf $\GG=\GG'(l)$ is generated by global sections, its  first three Chern 
classes are positive, while the  forth and following  ones with even indexes 
are negative. \\
  The sheaf $\FF=\GG\oplus \OPN^{n-2}$  is reflexive too, it is generated by 
global sections and it has the same Chern classes as $\GG$, so the \lq\lq 
even\rq\rq \ ones are negative from the fourth on.

  For the \lq\lq odd\rq\rq \ Chern classes, we can start from the rank 3 
reflexive sheaf $\GG'_1=\GG \oplus \OPN (a)$ (where $\GG$ is as above and  $a\gg 0$) such that 
$c_4(\GG'_1)>0$. Again, if  $l\gg 0$,    $\GG_1=\GG'_1(l)$ and 
$\FF_1=\GG_1 \oplus \OPN^{n-3}$  are generated by global sections, their first 
four Chern classes are positive, while the   fifth and following  \lq\lq 
odd\rq\rq \ ones  are negative.
\end{example}


\begin{thebibliography}*

\bibitem{fulton}  Fulton, W.: Intersection theory. Springer-Verlag, Berlin-Heidelberg (1984) 

\bibitem{alggeom} Hartshorne, R.: Algebraic Geometry. GTM 52, Springer-Verlag, New York (1977) 

\bibitem{HSRS}      Hartshorne, R.: Stable reflexive sheaves. 
Mathematische Annalen 254, 121-176 (1980) 

\bibitem{HSVB}  Hartshorne, R.: Stable vector bundles of rank 2 on 
$\PP^3$. Mathematische Annalen 238, 229-280 (1978) 

\bibitem{HL}    Huybrechts, D.,  Lehn, M.:  The geometry of moduli spaces of sheaves. Aspects Math. E31,  Bonn (1997). 

\bibitem{KL}  Kleiman, S. L.: The transversality of a general translate. Compositio Mathematica, 28 no. 3, p. 287-297, (1974).
 
\bibitem{MDP}  Martin-Deschamps, M., Perrin, D.: Sur la classification des courbes gauches, Ast\'erisque 184-185, (1990).

\bibitem{Mat}  Matsumura H.:  Commutative Algebra, Mathematics Lecture Notes Series, 56 ( 1980).


\bibitem{okted} Okonek, C.: Reflexive Garben auf $\PP^4$. Mathematische Annalen 260, 211-237 (1982). 

\end{thebibliography}
\end{document}